

\input amstex.tex
\input amsppt.sty

\catcode`\@=11
\def\logo@{}
\catcode`\@=\active

\parindent = 20 pt
\hsize = 6.5 true in
\vsize = 9   true in
\xspaceskip          = 7 pt plus 3 pt minus 0.2 pt


\define\uroman#1{%
        \uppercase
        \expandafter{\romannumeral#1}
        }

\font\bigrm = cmbx10 scaled \magstep 1
\define\ooi{{\dfrac 1 i}}
\define\tpi{{2 \pi i}}
\define\ztpi{\Z\cdot\tpi}
\define\ztpitau{\Z\cdot\tpi\tau}
\define\ootpi{\dfrac{1}{\tpi}}
\define\isom{\cong}

\def\O{{\Cal O}}
\define\Q{{\Bbb Q}}
\define\C{{\Bbb C}}

\define\R{{\Bbb R}}
\define\Z{{\Bbb Z}}
\def\P{{\Bbb P}}

\define\Sym{\operatorname{Sym}\nolimits}
\def\Re{\operatorname{Re}\nolimits}
\define\sgn{\operatorname{sgn}\nolimits}


\def\Im{\operatorname{Im}\nolimits}

\define\Gal{\operatorname{Gal}\nolimits}

\hyphenation{Spring-er}
\hyphenation{Ver-lag}
\hyphenation{theorem}

\define\divides{\mid}
\define\sumprime{\mathop{{\sum}'}}

\define\map#1{%
        \smash{\mathop{\longrightarrow}\limits^{#1}}}
\define\mapp#1#2{%
        \smash{\mathop{\longrightarrow}\limits^{#1}_{#2}}}
\define\mapdown#1{%
        \Big\downarrow
        \rlap{$\vcenter{\hbox{$\scriptstyle #1$}}$}}

\define\[{\left[}
\define\]{\right]}
\define\({\left(}
\define\){\right)}

\define\Wedge{\bigwedge}
\define\Dzeroq{D_{0,q}}
\define\realper{\Omega^+}
\define\sltwo{Sl_2(\Z)}
\define\gp{\Gamma_1(N)}
\define\zn{\zeta_N}
\def\bmod{\mathrel{\,\text{mod}\,}}
\define\Gab{G_{a,b}}
\define\sumab{\sum\limits_{
        \scriptstyle a \bmod N
        \atop
        \scriptstyle b \bmod N
        }}
\define\suma{\sum\limits_{  a \bmod N }}
\define\sumb{\sum\limits_{  b \bmod N }}
\define\cal{\Cal}
\def\L{{\cal L}}
\define\EN{{\cal E}_N}
\define\LN{{\cal L}_N}

\define\lattice#1#2{\Z \cdot #1 + \Z \cdot #2}
\define\leftquot{\setminus}
\define\proof{\demo{\bf Proof}}
\def\qed{{\bf Q.E.D. }}
\define\endproof{\qed\enddemo\smallskip}
\define\half{{1 \over 2}}

\define\hdg#1#2{                        %
         \line                          %
          {                             %
          #1.                           
          \hskip 1 em                   %
          #2                            
          \hss                          %
          }                             %
        }

\define\section#1.#2\par{                       %
        \vskip 24pt plus 5 pt minus 2 pt        %
        \penalty -5000                          %
        {                                       %
         \bf
         \hdg{#1}{#2}                           %
        }                                       %
        \bigskip
        \message{Section #1}%
}

\define\relunder#1\under#2{\mathrel{\mathop{\kern 0pt \relax #2}\limits_{#1}}}
\def\goodbreak{\penalty-1000}
\define\today{\number\year\ \ifcase\month\or January\or February\or March\or April\or
        May\or June\or July\or August\or September\or October\or
        November\or December\fi\ \number\day}


\bigskip

\centerline {\bigrm Eisenstein series of weight one, $q$-averages of the $0$-logarithm and periods
   of elliptic curves}
\centerline {by}
\centerline {%
        Daniel R. Grayson%
        \footnote{%
                University of Illinois at Urbana-Champaign.
                Research supported by the NSF.
                }
        and Dinakar Ramakrishnan%
        \footnote{%
                California Institute of Technology.
                Research supported by the NSF.
                }
        }
\bigskip
\rightline{To Kumar Murty:  Hyapi aravai}
\bigskip

{\bf Abstract}

\medskip

For any elliptic curve $E$ over $k\subset \R$ with $E(\C)=\C^\times/q^{\Z}$, $q=e^{2\pi iz}, \Im(z)>0$, we study the $q$-average $D_{0,q}$, defined on $E(\C)$, of the function $D_0(z) = \Im(z/(1-z))$. Let $\Omega^+(E)$ denote the real period of $E$. We show that there is a rational function $R \in \Q(X_1(N))$ such that for any non-cuspidal real point $s\in X_1(N)$ (which defines an elliptic curve $E(s)$ over $\R$ together with a point $P(s)$ of order $N$), $\pi D_{0,q}(P(s))$ equals $\Omega^+(E(s))R(s)$. In particular, if $s$ is $\Q$-rational point of $X_1(N)$, a rare occurrence according to Mazur, $R(s)$ is a rational number.

\section 1. Introduction

The relationship between modular forms of weight one and periods of elliptic
curves is well-known, certainly to experts in the field.  In this paper
we study in some detail
those modular forms of weight one which arise from $q$-averages of the
$0$-logarithm
function, $\ell_0(z) = {z/(1-z)}$.
We are led to consider these rather special forms to point out
the analogy with the situation in algebraic $K$-theory, where
the other polylogarithm functions have already played an important role in
connection with special values of $L$-functions.

Let E be an elliptic curve over $\Q$ and let $P_0 \in E(\Q)$ be a rational
point of order $N>1$.  Writing $E(\C) = \C/L$ for some complex lattice $L =
\Z\cdot1 + \Z\cdot\tau$, with $\Im \tau > 0$, we let $q = e^{\tpi \tau }$;
the exponential map $ u \mapsto e^{\tpi u}$ gives an isomorphism $\Phi :
E(\C) \isom \C^\times/q^\Z$.  Let $z_0$ be a representative for the coset
$\Phi(P_0)$, so that $\Phi(P_0) = z_0\cdot q^\Z$.

For $k \ge 0$ define the polylogarithm functions,  $$\ell_k(z) =
\sum_{n=1}^\infty{z^n \over n^k};$$
these functions can be regarded as single-valued functions on the disk $\{z : |z|
< 1 \}$, or they can be regarded, by analytic continuation, as multi-valued
functions on $\{ z : z \ne 1, \, z \ne 0 \}$.
The Bloch-Wigner dilogarithm function
$$
D_2(z)
= \log|z| \arg(1-z) + \Im \ell_2(z)
= \log|z| \arg(1-z) - \Im \int_0^z \log(1-t){dt \over t}
$$
is a single-valued real function of the complex
variable $z$.  Summing over a coset of $q^\Z$ yields $D_{2,q}(z) =
\sum_{n \in \Z} D_2(z q^n)$, which is a real function on $E(\C)$.  Its value at
$P_0$ plays an important role in conjectures of Spencer Bloch about values of
the L-series $L(E,s)$ at $s=2$ in connection with values of the higher
regulator map $K_2E \rightarrow \R$ from algebraic K-theory (\cite{4}, \cite{2}) and certain
Eisenstein-Kronecker-Lerch series (\cite{26}).
In particular, it often happens that $${ \pi D_{2,q}(z_0)
\over L(E,2) } \tag1.1$$ is a rational number to great numerical accuracy.

In this paper we explore the analogous situation at $s = 1$.

To do this,
we replace the dilogarithm by
the first member of the sequence of polylogarithm functions, namely, the
$0$-logarithm.  We let $$D_0(z) = \Im \ell_0(z) = \Im {( z /( 1-z ))}
= \Im { (1/( 1-z ))}$$ for $z
\not= 1$, and let $D_0(1) = 0$.  Summing over a coset gives a function
$$\Dzeroq(z) = \sum_{n \in \Z} D_0(z q^n) $$ which amounts to a real function on
$E(\C)$.  Something similar works for the higher polylogarithm functions
$\ell_m(z)$, but other
terms must be added to arrive at a function which is single-valued and
real-analytic on the complex plane, generalizing the $m=2$ case of \cite{3},
cf. \cite{21}, \cite{28}, where the latter presents Zagier's conjecture relating to the values of Dedekind zeta functions of arbitrary number fields $F$ at positive integers $m\geq 2$ and $K_{2m-1}(F)$.

The conjecture of Birch and Swinnerton-Dyer \cite{24, p.~362}
states, among other things, that the coefficient of the leading term in the
Taylor series of $L(E,s)$ at $s=1$ is a rational number times the product of
the real period of $E$ and the determinant of the height pairing on the
Mordell-Weil group $E(\Q)$, with the rank of the Mordell-Weil group and the order
of vanishing of $L(E,s)$ at $s=1$ predicted to be the same. In particular,
$L(E,1)$ should always be a
rational multiple (possibly zero) of the real period of $E$, which has long been known in
the case where $E$ is a modular elliptic curve (via the theory of
modular symbols and the Manin--Drinfeld theorem \cite{9}).
Now every elliptic curve over $\Q$ is known to be modular by Wiles, Taylor,
Diamond, Conrad and Brueil (\cite{5}, \cite{25}, \cite{27}).
Thus, in the search for an analogue of the rationality
of (1.1), we may as well use the real period instead of $L(E,1)$;
this even simplifies the analysis, for ignoring the arithmetic data that underlies
$L(E,1)$ allows us to study the variation as $E$ moves in a family. (It is not important to what follows, but it may be worthwhile to remark that by Kolyvagin \cite{16}, $E(\Q)$ is finite when $L(E,1)=0$, and one of the ingredients of his theorem was a non-vanishing hypothesis, established later by two sets of authors, one of them being Kumar and Ram Murty, and the other being Bump, Friedberg and Hoffstein.)

We write the
equation for E as $y^2 = f(x)$ with $f$ a monic polynomial of degree 3.
Let $\omega = {dx /( 2 y)}$ be the standard
holomorphic differential on $E(\C)$; it depends on the choice of equation.
Set the real
period to be $$\realper = \int_{E(\R)^\circ} \omega = \int_\gamma^\infty{ dx
\over \sqrt{f(x)}} $$
where $\gamma$ is the largest real root of $f(x)$.

We introduce the group
$$
\gp = \biggl\{ \pmatrix a&b\cr c&d\endpmatrix \in \sltwo
 \biggm|
    \pmatrix a&b\cr c&d\endpmatrix \equiv \pmatrix 1&*\cr 0&1\endpmatrix \mod N
 \biggr\}.
$$
Let $X_1(N)$ denote the modular curve corresponding to the group $\gp$.  It is
a projective curve
defined over $\Q$, so that $X_1(N)\otimes\C$ contains the quotient
$\gp\leftquot H$ as an
affine open subset.  For each $s \in X_1(N)$, not a cusp, there is a
corresponding elliptic curve $E(s)$ equipped with a point $P(s)$ of order $N$;
the curve and point are both defined over $\Q(s)$.
It is possible to pick a differential $\omega(s) \ne 0$ on $E(s)$ that
varies algebraically in $s$; this enables us to interpret
the various quantities
above, such as $q$ and $\Omega^+$, as functions of~$s$.

\proclaim{Theorem 1.2} Let $N > 0$.  There is a
rational function $R \in \Q(X_1(N))$ such
that for each real point $s \in X_1(N)$, not a cusp, we have $R(s) = 2 \pi
D_{0,q(s)}(P(s)) / \Omega^+(s)$.
\endproclaim

In particular, if $s$ is a rational point, then the quantity $R(s)$ is a
rational number.

This theorem was suggested by two sequences of numerical experiments.
Our first sequence of experiments
showed that $$ { \pi \Dzeroq(P_0) \over \realper} $$
appears to be a rational number of small height
to great accuracy, for example, 100 digits.

\topinsert
$$
\vcenter{%
  \offinterlineskip
  \def\tail{&}
  \def\head{&}
  \def\bitofspace{height 3 pt&\omit&&\omit&&\omit&\cr}
  \def\\{\tail\cr\bitofspace\noalign{\hrule}\bitofspace\head}
  \def\CR{\tail\cr\head}
  \hrule
  \tabskip=0pt
  \halign to 4 in {%
    \vrule#\tabskip = 1em plus 2em&%
    \vrule height 9.5pt depth 4.5pt width 0pt
    \hfil$#$\hfil&%
    \vrule#&%
    \hfil$#$\hfil&%
    \vrule#&%
    \hfil$#$\hfil&%
    \hfil#\vrule\tabskip = 0pt\cr
   \bitofspace
   \head
    N && k && 2 N R = 4 N \pi \Dzeroq(k P) / \realper                   \\
    3 && 1 &&  -a_1 \rlap{\hbox{ \quad (with $a_1 = a_3$)}}             \\
    4 && 1 &&  -2                                                       \\
    5 && 1 &&  c - 3                                                    \CR
      && 2 &&  - 3 c - 1                                                \\
    6 && 1 &&  -4                                                       \\
    7 && 1 &&  d^2 -3 d - 3                                             \CR
      && 2 &&  -5 d^2 + d + 1                                           \CR
      && 3 &&  3 d^2 - 9 d + 5                                          \\
    8 && 1 &&  (-8 d + 2)/d                                             \CR
      && 3 &&  (-8 d + 6)/d                                             \\
    9 && 1 &&  f^3-3f^2-5                                               \CR
      && 2 &&  -7f^3+3f^2-1                                             \CR
      && 4 &&  -5f^3+15f^2-18f+7                                        \\
   10 && 1 &&  - 4f - 4                                                 \CR
      && 3 &&  -12f + 8                                                 \\
   12 && 1 &&  { (24 \tau^3 - 48 \tau^2+36\tau-10)/(1-\tau)^3}          \CR
      && 5 &&  {(24\tau^3-24\tau^2+12\tau-2)/(1-\tau)^3}                \tail\cr
      \bitofspace\noalign{\hrule}
    }
  }
$$
\centerline{\bf Table 1.}
\endinsert

Our second sequence of numerical experiments was this.  The possible numbers N that occur
are $ 1 \leq N \leq 10 $ or $ N = 12 $
(according to a theorem of Mazur, \cite{19}), and for
each of these there is a one-parameter family (defined over $\Q$) of pairs
$(E_t,P_t)$ where $P_t$ is a point of order N on the elliptic curve $E_t$.  We
found that the ratio above is a simple function of the modular parameter $t \in
\C$.
The results are summarized in Table~1, which presents, for each case, a
rational function which is a good fit to numerical results with many digits of
accuracy, for a finite number of values of $t$.
Refer to \cite{17, p.217}, for the parametrizations used; we use
Kubert's names
for the parameters.  In each case, the family of elliptic curves is represented
by an equation of the form $y^2 + a_1 x y + a_3 y = x^3 + a_2 x^2 + a_4 x +
a_6$; we choose $\omega = dy/(2y + a_1 x + a_3)$, as is customary.

We expected that the ratio $R$ would be a rational function of the parameter,
i.e., would be a modular function.  We were surprised that its coefficients
were always rational, but we shouldn't have been, for it turns out that the
well-known Hecke-Eisenstein modular forms of weight one have sufficient
rationality properties to explain the observations.  These forms are presented
in Lang's book \cite{18} in chapter XV, section 1, in a
concise form.  The best
explanation of their rationality properties is presented by Nicholas Katz in
\cite{14, Appendix C}, where he presents a purely algebraic
construction of the Hecke-Eisenstein modular forms; we hadn't been aware of his work during our initial investigations.
His construction shows the forms are defined
over the rational numbers, and are computable exactly
and algebraically, so we were able to verify that the entries in Table 1 are
correct by means of a second computer program.

In this paper we explain how to deduce the rationality of the ratios $R$
presented above from the rationality of the Hecke-Eisenstein forms.  This is
done by using $q$-expansions to identify our ratio with the modular
form presented in Katz's paper.
This identification is easy in the case where $P_0$ is on the identity
component of $E(\R)$, but harder otherwise, essentially because Katz provides
the $q$-expansion just at certain cusps.
  We then explain the relationship with the
Hecke-Eisenstein forms as presented in Lang's book.  We present an exposition
of the $q$-expansion principle which is sufficient to deduce the desired
rationality of the ratios $R$ directly from the rationality of the
$q$-expansion coefficients, allowing one to bypass, if desired, the elegant
algebraic construction of Katz.

We end with a general question.  Let $M$ be a motive over $\Q$, of weight $k\ge
1$, defined by an irreducible direct factor of the cohomology of an abelian
variety $A$.  If $r$ is an integer with $k/2 \le r \le k+1$, critical or not,
is there an analog of the function $D_{0,q}$ (or $D_{2,q}$) on $A(\C)$, whose
values at torsion points have (at times) rational relations with $L(M,k)$?  For
$M = \Sym^2(E)$, with $E$ an elliptic curve over $\Q$, k=2, and $r=3$, one has
such a phenomenon; see the papers \cite{8} and
\cite{20}.

We thank the referee for helpful comments on an earlier version, which have been incorporated.

\section 2. Interpreting $\Dzeroq$ as a modular form of weight one for $\gp$


We use the same notation as in the introduction.
Since $E$ is real, we may choose the isomorphism $E \cong \C/L$
in such a way that the real structure on $E$ corresponds to complex conjugation
on $\C$.  Then since $L$ is invariant under complex conjugation we may choose
$\tau$ so that $\Re \tau$ is either $0$ or $1/2$; in the first case, $E(\R)$
has two components and $q > 0$, and in the second case, $E(\R)$ has one
component and $q < 0$.  In both cases $|q| < 1$.
Since $P_0$ is a point of order $N$, we may write $z_0$ in the form $z_0 =
\zn^k q^{\ell/N}$ where $\zn = \exp(\tpi / N)$.
Since $P_0$ is real, we may pick $z_0$ so that either

\medskip

        case A) $z_0 = \zn^k$

\medskip

\noindent or

\medskip

        case B) $z_0 =  \zn^k q^{1/2} $.

\medskip

\noindent  Case B occurs when
$2 \divides N$, $q > 0$, and $P_0$ is not
on the identity component of $E(\R)$.

        We occupy ourselves with case A first.  Making use of the easy
identities $D_0(\bar z) = -D_0(z) = D_0(1/z)$ and remembering that q is real
we compute

$$
\eqalign{
  \ooi \Dzeroq(\zeta_N^k)
     &=  \ooi \sum_{j\in \Z} D_0(\zn^k q^j)
  \cr& = \ooi D_0(\zn^k) + { 2 \over i }  \sum_{j=1}^\infty D_0(\zn^k q^j)
  \cr&= \ooi D_0(\zn^k) + {2\over i} \Im \sum_{j=1}^\infty \sum_{m=1}^\infty (\zn^k q^j)^m
  \cr& = \ooi D_0(\zn^k) + {2\over i} \Im \sum_{n=1}^\infty \biggl(
        \sum_{\scriptstyle m\divides n \atop \scriptstyle m>0 } \zn^{k m}
        \biggr) q^n
  \cr&= \ooi D_0(\zn^k) - \sum_{n=1}^\infty \biggl(
        \sum_{\scriptstyle m\divides n \atop \scriptstyle m>0 }
                 (\zn^{k m} - \zn^{-k m})
        \biggr) q^n
\cr}
\tag 2.1
$$
The double sums above are absolutely convergent because $|q| < 1$, and
$$
        \sum_{j=1}^\infty \sum_{m=1}^\infty |q|^{j m} \le { |q| \over (1-|q|)^2.
}
$$

We set
$$
g_k(\tau) = \ooi D_0(\zn^k) - \sum_{n=1}^\infty \biggl(
        \sum_{\scriptstyle m\divides n \atop \scriptstyle m>0 }
                 (\zn^{k m} - \zn^{-k m})
        \biggr) q^n
$$
so that $g_k(\tau) = (1/i) \Dzeroq(\zn^k)$ when $q$ is real.
We will use this $q$-expansion for values of $q$ which are not necessarily
real, as $\Dzeroq(z_0)$ is inappropriate when $q$ is not real.

Now we consider case B.  We have
$$
\eqalign{
        \ooi \Dzeroq(\zn^k q^{1/2})
        &= \ooi \Dzeroq\left( \zn^k q^{{1 \over 2}} \right)
         = \ooi \sum_{j \in \Z} D_0\left( \zn^k q^{{1 \over 2}+j} \right)
        \cr&= {2\over i} \sum_{j=1}^\infty D_0\left( \zn^k q^{j-{1 \over 2}}\right)
            = {2 \over i} \sum_{j=1}^\infty \Im
                \sum_{m=1}^\infty \zn^{k m} q^{(j-{1 \over 2})m}
        \cr&= {2 \over i} \sum_{n=1}^\infty \biggl( \sum_{
                \scriptstyle m\divides n
                \atop
                {
                \scriptstyle n/m \text{\hskip 3pt odd}
                \atop
                \scriptstyle m > 0
                }
                } \Im \zn^{k m} \biggr) q^{n/2}
        \cr&= \sum_{n = 1}^\infty \Biggl(
                {2\over i} \sum_{
                \scriptstyle m\divides n
                \atop
                {
                \scriptstyle m \text{\hskip 3pt odd}
                \atop
                \scriptstyle m > 0
                }
                } \left( \Im \zn^{k n / m}\right) \Biggr)
                q^{n/2}
        \cr&= \sum_{n = 1}^\infty \Biggl(
                - \sum_{
                \scriptstyle m\divides n
                \atop
                {
                \scriptstyle m \text{\hskip 3pt odd}
                \atop
                \scriptstyle m > 0
                }
                } \left( \zn^{k n / m} - \zn^{-k n / m}  \right) \Biggr)
                q^{n/2}.
\cr}
\tag 2.2
$$

We set
$$
h_k(\tau) =
 \sum_{n = 1}^\infty \Biggl(
                - \sum_{
                \scriptstyle m\divides n
                \atop
                {
                \scriptstyle m \text{\hskip 3pt odd}
                \atop
                \scriptstyle m > 0
                }
                } \left( \zn^{k n / m} - \zn^{-k n / m}  \right) \Biggr)
                q^{n/2}
$$
so that $h_k(\tau) = (1/i) \Dzeroq(\zn^k q^{1/2}) $ when $q$ is real.

        Let $\L$ denote the set of all lattices $L$ in $\C$.
Let $\LN$ denote the set of all pairs $(L,u)$ where $L$ is a lattice in
$\C$ and $u$ is an element of $N^{-1}L$ of order $N$ in $N^{-1}L/L$
(see \cite{18, p.~101}).
Consider a function $F :
\LN \rightarrow \C$ satisfying the identity $F(\alpha L, \alpha u) =
\alpha^{-\ell} F(L,u)$ for all $\alpha \in \C$.
Call such a function {\sl homogeneous of degree $-\ell$}.
Since k is relatively prime to
N, we have a map $\phi_k : H \rightarrow \LN$ from the upper half plane $H$,
defined by
$$
\phi_k(\tau) = (\Z\cdot \tpi + \Z \cdot \tpi \tau , \tpi k/N).
$$
Notice the factor of $\tpi$ used here.
The composite $ f_k =
F \circ \phi_k$ is a modular form of weight $\ell$ for the group $\gp$
if it is meromorphic
on the upper half plane and
 at the cusps; the main import of being a modular form is
the identity
$$
f_k({a \tau + b \over c \tau + d}) = (c \tau + d)^{\ell} f_k(\tau)
$$
for all $\pmatrix a&b\cr c&d\endpmatrix\in \gp$, which follows from the
homogeneity of $F$.

For $f_k$ to be meromorphic at the cusps means that
$$
(c\tau + d)^{-\ell} f_k\({ a \tau + b \over c \tau + d }\)
$$
is a meromorphic function of $e^{\tpi\tau/N}$, for each $\pmatrix a&b\\c&d
\endpmatrix \in \sltwo $
(this condition is formulated incorrectly in \cite{18, p.~103}).  This
condition is independent of $k$, and in terms
of $F$ it means that
$$
F\(\ztpi+\ztpitau, \tpi\({c\tau+d\over N}\)\)
$$
is a meromorphic function of $e^{\tpi\tau/N}$ for each pair of relatively
prime intergers $c$, $d$.

The field of definition of a modular form is properly understood as described
in \cite{14}.
Given a scheme $S$ on which the integer $N$ is invertible,
we let $\EN(S)$ denote the set of triples $(E,\omega,P)$, where $E$ is a family
of elliptic curves over $S$, $\omega$ is a holomorphic differential form on $E$
relative to $S$ which is nonzero on each fiber, and $P$ is a section of $E$
over $S$ which has order $N$.  There is a bijection $\LN \to \EN(\C)$ which
sends $(L,u)$ to $(\C/L,dz,u+L)$, where $z$ is the coordinate function on $\C$.

If $R$ is a ring, then a holomorphic modular form for $\gp$ of weight $\ell$
over $R$ is a collection of maps $F : \EN(S) \to H^0(S,\O_S)$ defined for any
$R$-scheme $S$, which is natural in $S$, and which is homogeneous of weight
$\ell$ in the sense that $F(E,\alpha\omega,P) = \alpha^{-\ell}F(E,\omega,P)$
for each $\alpha \in H^0(S,\O_s)$.  When $R$ is a field, then we say that $R$
is a field of definition of $F$.

When $\C$ is an $R$-algebra, we will make use of the bijection mentioned
previously implicitly, regarding such an $F$ also as a function $F : \EN(\C)
\to \C$, and writing $F(\C/L,dz,u+L) = F(L,u)$.

The $q$-expansions at the cusps are obtained by choosing $k$ and $\ell$ so
$\gcd(k,\ell,N) = 1$ and expanding
$$
F\(\ztpi+\ztpitau, \tpi\({k+\ell\tau\over N}\)\)
$$
as a Laurent series in $q^{d/N}$, where $d = \gcd(\ell,N)$.

In \cite{14, p.~260} is described a modular form, $A_1$, of weight
$1$.  Some of its
$q$-expansions are given in \cite{14, (2.7.10)} as
follows.
$$
        A_1(\Z \cdot \tpi + \Z \cdot \tpi \tau, \tpi k / N) =
        {1 \over 2 i} \cot(\pi k / N)
        - \sum_{ n=1}^\infty \biggl(
                \sum_{\scriptstyle m\divides n \atop \scriptstyle m>0 }
                         (\zn^{k m} - \zn^{-k m})
                \biggr) q^n
\tag 2.3
$$
One sees easily that
$$
        {1 \over i} D_0(\zn^k)
 =
        {1 \over i} \Im \( { \zn^k \over 1 - \zn^k } \)
 =
        {1 \over 2} \( { \zn^k + 1 \over \zn^k - 1 } \)
 =
        {1 \over 2 i} \cot(\pi k / N).
$$
Thus from the identity of the $q$-expansions (2.1) and (2.3),
we see that for case A
$$
g_k(\tau) =
A_1(\Z\cdot \tpi + \Z \cdot \tpi \tau , \tpi k/N).
\tag 2.4
$$

In case B, we expect the analogous identity
$$
        h_k(\tau)
 =
        A_1\bigl(\Z\cdot \tpi + \Z \cdot \tpi \tau ,
                \tpi\bigl( {k\over N} + \half \tau \bigr)\bigr).
\tag 2.5
$$
to hold
(it doesn't seem to follow immediately from \cite{14, 2.7.8});
proving it amounts to computing the $q$-expansion for $A_1$ at cusps
other than the ones Katz considered, and we do this later in (3.6).

In \cite{14, Appendix~C} is an amazing algebraic construction of
$A_1$ which shows that it is defined over $\Z[{1 \over 6 N}]$.  We will use
this result to get the rationality.

Let $E$ be an elliptic curve, as in the introduction, defined over a subfield
$K$ of $\C$, with a point $P_0$ of order $N$ and a differential $\omega =
dx/(2y)$.  We can find a lattice $L$, a number $u_0 \in \C$, and an isomorphism
$(E,\omega,P_0) \isom (\C/L,dz,u_0+L)$.  Since $A_1$ is defined over $\Q$, we
see that $r := A_1(L,u_0) \in K$.  Now write $L = \Z \cdot v_1 + \Z \cdot v_2$ with
$\Im\(v_2/v_1\) > 0$, and let $\tau := v_2/v_1$.  Thus
$$
\eqalign{
{ \tpi \over v_1} & A_1(\Z\cdot \tpi + \Z \cdot \tpi \tau, \tpi u_0/v_1)
=
{1 \over v_1} A_1(\Z + \Z \cdot \tau, u_0/v_1)
\cr&=
A_1(\Z\cdot v_1 + \Z \cdot v_2, u_0)
=
A_1(L,u_0)
=
r
\in K.
\cr
}
$$

Now suppose we are in case A, so that $K \subseteq \R$,
we take $\realper = v_1$ to be real and positive, and we pick $k \in \Z$ so
$k/N = u_0/v_1$.  Then from (2.4) we deduce:
$$
\eqalign{
\cr
{2 \pi \Dzeroq(\zn^k) \over \realper }
&= \( \ooi \Dzeroq(\zn^k) \) \Bigm/ \( \ootpi \realper \)
 = \( g_k(\tau) \) \Bigm/ \( \ootpi \realper \)
\cr
&= { \tpi \over v_1} A_1(\Z\cdot \tpi + \Z \cdot \tpi \tau, \tpi k/N)
= r \in K
\cr}
$$
This shows the desired rationality of the values of $R$ presented in the
introduction.

In case B, we have ${k \over N} + {1 \over 2} \tau = { u_0/ v_1 }$, so by
(2.5) we have:
$$
\eqalign{
{2 \pi \Dzeroq(\zn^k q^{1/2}) \over \realper }
    &= \( \ooi \Dzeroq(\zn^k q^{1/2}) \) \Bigm/ \( \ootpi \realper \)
\cr &= h_k(\tau) \Bigm/ \( \ootpi \realper \)
\cr &= {\tpi\over v_1}A_1(\ztpi+\ztpitau,\tpi(k/N+\tau/2))
= r \in K,
\cr}
$$
yielding the desired rationality in this case.

The sense in which $R$ is a rational function of the parameter is this.  Each
of the families of elliptic curves in Table~1 can each be viewed as an element
$(E,\omega,P_0) \in \EN(S)$, where $S$ is an open subset of $\P^1_\Q$.  Thus $R
= A_1(E,\omega,P_0) \in \Q(\P^1)$ is a rational function.  Alternatively, for
$N \ge 3$, one may take $S = X_1(N) - \{\text{cusps}\}$ and $(E,\omega,P)$ the
universal family of elliptic curves, thereby proving theorem
1.2.

\section 3. Hecke-Eisenstein modular forms

Now we recall work of Hecke (\cite{10}, \cite{11}).
(It also occurs in \cite{12, Chapter 3}.)
We always use an unadorned $\equiv$ to denote
congruence modulo N.

For $(L,u) \in \LN$ and $s \in \C$ with $\Re s > 1 $ we define
$$
\Phi(L,u,s) = \mathop{{\sum}'}_{{\omega / N} \equiv u \mod L} \omega^{-1}
|\omega|^{-s},
$$
where the prime means that the term with $\omega=0$ is to be omitted if it
occurs.
We define
$$
G(L,u) = \Phi(L,u,0)
$$
by analytic continuation as Hecke does.  It is evident that $G$ is a
homogeneous
function of degree $-1$.
For $a,b \in \Z$ we let $G_{a,b}(\tau) = G\(\Z\tau+\Z,(a \tau + b)/N\)$,
in
accordance with the notation in \cite{22}.  If it happens that $a=0$
then $G_{a,b}(\tau)$ is a modular function for $\gp$, but in any case, it is a
modular function for $\Gamma(N)$.  Providing a $q$-expansion at $\tau = \infty$
for each function $G_{a,b}(\tau)$ in terms of $e^{\tpi\tau/N}$ is the same as
providing a $q$-expansion for $G$ at each of the cusps of $X(N)$.  Conversely,
finding a $q$-expansion for $G_{a,b}(\tau)$ at $\tau=\infty$ for all $a,b$ is
the same as finding a $q$-expansion for one of the functions $G_{a,b}(\tau)$ at
all of the cusps.
Now we present some notation needed to write down the $q$-expansion found by
Hecke.

Introduce the Hurwitz zeta
function
$$
\zeta(s,\alpha) = \sum_{n>-\alpha} (n+\alpha)^{-s}
$$
and the notation
$$
\delta(v) = \cases
1,    &   \text{if $v \in \Z$; } \cr
0,    &   \text{otherwise.}
\endcases
$$
 Define
$$
\alpha_ n(a,b) =
\cases
        \displaystyle
        {1 \over N} \delta( { a \over N } )
                \lim_{s\rightarrow 1} \left[
                        \zeta(s,{b \over N}) - \zeta(s,{-b \over N})
                        \right]
                - {\pi i \over N} \left[
                        \zeta(0,{a \over N}) - \zeta(0,{-a \over N})
                        \right]
        & \text{\quad if $ n = 0$}
\cr
        & \strut
\cr
        \displaystyle
        - { \tpi \over N } \sum_{
                \scriptstyle  m\divides n
                \atop
                \scriptstyle{ \scriptstyle  n \over \scriptstyle m } \equiv a
                }
            (\sgn m) \zn^{b m}
        & \text{\quad if $ n > 0$}
\endcases
$$

\bigskip\goodbreak

Here is Hecke's result.

\proclaim{Proposition 3.1} {\rm \cite{22, p.~168}, or
\cite{11, p.~203}.}
$
G_{a,b}(\tau) = \sum\limits_{ n = 0}^\infty \alpha_ n(a,b) q^{ n / N}
$
\endproclaim

        From this we can deduce the following.

\proclaim{Proposition 3.2}
$
\ootpi \sumab \zn^{k a} G_{a,b}(\tau) = g_k(\tau)
$
\endproclaim

\proof
Part of this proof is just like the proof of \cite{14, (2.7.12)}.

We expand the left-hand side using
Proposition 3.1, obtaining
$$\sum\limits _{ n = 0} ^\infty\beta_ n q ^{ n/N},$$
where
$$\beta_ n = {1 \over \tpi} \sumab \zn^{k a} \alpha_ n(a,b).$$

        We take the case $ n = 0$ first.  We have
$$
\eqalign{
        \beta_0
        &= {1 \over \tpi}
           \sumab \zn^{k a} \alpha_0(a,b)
        \cr&= - {1 \over 2} \suma \zn^{k a} \left[
           \zeta(0,{a \over N}) - \zeta(0,{-a \over N}) \right]
        \cr& = - {1 \over 2} \suma ( \zn^{k a} - \zn^{- k a} ) \zeta(0,{a \over N})
\cr}
$$
We see from \cite{26, p.59}, that $\zeta(0,\alpha) = \half - \alpha$ for
$0<\alpha<1$, and it is known that $\zeta(0,1) = \zeta(0) = -\half = \half -
1$.  Thus (compare with \cite{14, p.~263}) we have
$$
\eqalign{
           \beta_0
        &= - {1 \over 2} \sum_{a=1}^N (\zn^{k a} - \zn^{-k a})(\half - {a \over N})
             = {1 \over 2 N} \sum\limits_{a=1}^N (\zn^{k a} - \zn^{-k a})a
        \cr& = - {1 \over N i} \Im \sum\limits_{a=1}^N a \zn^{k a}
             = \ooi \Im \lim_{T \to \zn^k}
                { \displaystyle \sum_{a=1}^N T^a \over \displaystyle 1 - T^N }
        \cr& = \ooi \Im \lim_{T \to \zn^k} { T \over 1-T }
             = \ooi \Im { \zn^k \over 1 - \zn^k }
             = \ooi D_0(\zn^k).
        \cr
}
$$
(The middle equality above is an application of l'H\^opital's rule.)
This agrees with the coefficient of $q^0$ in $g_k(\tau)$.

        Now we take up the case $ n > 0$.  We have
$$
\eqalign{
  \beta_ n
    &= - {1 \over N} \sumab \zn^{k a} \sum_{
                \scriptstyle  m\divides n
                \atop
                \scriptstyle{ \scriptstyle  n \over \scriptstyle m } \equiv a
                }
        (\sgn m) \zn^{b m}
    \cr& = - { 1 \over N} \sum_{m\divides n} \zn^{k  n / m}(\sgn m) \sumb \zn^{b m}
    \cr&= - \sum _{
                \scriptstyle m\divides n
                \atop
                \scriptstyle m \equiv 0
                } (\sgn m) \zn^{k  n / m}
\cr}
$$
and so $\beta_ n = 0$ unless $ n \equiv 0$, for else the summation has no terms.  We compute, for $ n > 0$,
$$
\eqalign{
        \beta_{N  n}
        &= - \sum_{m\divides n} (\sgn n) \zn^{k  n / m}
         = -  \sum_{m\divides n} ( \sgn m) \zn^{k m}
        \cr&= - \sum_{
                        \scriptstyle m \divides  n
                        \atop
                        \scriptstyle m > 0
                  } ( \zn^{k m} - \zn^{- k m})
\cr}
$$
This is the coefficient of $q^n$ in the expansion for $g_k(\tau)$. \endproof

        Now we show how to use Proposition 3.2 to
interpret $g_k(\tau)$
as a holomorphic modular form of weight one for the group $\gp$.

As standard basis for the exterior power $\Wedge^2_{\R}\C$ we will use $1\wedge i$.  We let
$\alpha_L$ denote the generator of the group $\Wedge^2_{\Z}L$ which is a
positive
multiple of $1\wedge i$.
For any $u,v \in N^{-1}L$ we have
$$
N{u \wedge v \over \alpha_L} = N^{-1} {Nu \wedge Nv \over \alpha_L}
        \in N^{-1}\Z
$$
and thus
$$
 \chi^L_u(v) \relunder \text{def} \under = \exp \left( \tpi N {
        u \wedge v
        \over
        \alpha_L
        } \right)
$$
is an $N^{th}$ root of 1.
We remark that $\chi^L_u(v) = e_N(u,v)$ in terms of the $e_N$-pairing (see
\cite{12, p.~477}).
One checks the
following formulas.
$$
\vcenter{
\openup 3\jot
\ialign{%
 \strut\hfil$\displaystyle{#}$&$\displaystyle{{}#}$\hfil&\qquad#\hfil\cr
 \chi^L_{u+u'}(v) &= \chi^L_u(v) \chi^L_{u'}(v)          &                      \cr
 \chi^{\alpha L}_{\alpha u}(\alpha v) &= \chi^L_u(v)     & for $\alpha \in \C$  \cr
 \chi^L_u(v+v')   &= \chi^L_u(v) \chi^L_u(v')            &                      \cr
 \chi^L_u(v)      &= 1                                   & if $u \in L$         \cr
 \chi^{\Z + \Z \cdot \tau}_{(k+\ell\tau)/N}\bigl({a \tau + b\over N}\bigr)
                  &= \zn^{ka-\ell b}                     & for $a,b,k \in \Z$   \cr
 }
}
$$
\noindent It follows that $\chi^L_u$ depends only on the class of $u$ in $N^{-1}L/L$, and
is a homomorphism $N^{-1}L \rightarrow \C^\times$.


For $(L,u) \in \LN$ we define
$$
G'(L,u) = \sum_{w\in N^{-1}L/L} \chi^L_u(w) G(L,w);
$$
the function $G'$ is a meromorphic modular form for $\gp$.  From
(3.2) we see that
$$
g_k(\tau) = G'(\ztpi+\ztpitau,\tpi k/N),
\tag 3.3
$$
which shows that $g_k(\tau)$ is a holomorphic modular form for $\gp$.
Combining this with (2.4) we see that
$$
A_1 = G'
\tag 3.4
$$

\goodbreak

        We can also deduce the following proposition.

\proclaim{Proposition 3.5} If $ 2 \divides N$, then
$
\ootpi \sumab \zn^{k a} (-1)^b G_{a,b}(\tau) = h_k(\tau)
$
\endproclaim

\proof
We expand the left hand side,  using
Proposition~3.1,
obtaining
$$
\sum_{ n=0}^\infty \gamma_ n q^{ n/N},
$$
where
$$
\gamma_ n = { 1 \over \tpi} \sumab \zn^{k a}(-1)^b\alpha_ n(a,b).
$$
We compute (keeping in mind that N is even)
$$
\eqalign{
        \gamma_0 &= {1 \over \tpi} \sumab \zn^{k a} (-1)^b \alpha_0(a,b)
\cr
        &= {1 \over \tpi}\biggl[ \sumb (-1)^b{1 \over N}
                \lim_{s\rightarrow 1} \left\{
                        \zeta(s,{b \over N}) - \zeta(s,{-b \over N})
                        \right\}
\cr
        &\quad - { \pi i \over N} \sumab \zn^{k a} (-1)^b \left\{
                         \zeta(0,{a \over N}) - \zeta(0,{-a \over N})
                         \right\}
                \biggr]
\cr
        &= 0
\cr}
$$
Now for $ n > 0$ we find
$$
\eqalign{
        \gamma_ n &= { -1 \over N} \sumab \zn^{ k a} (-1)^b
                \sum_{
                        \scriptstyle m\divides n
                        \atop
                        \scriptstyle { n \over m} \equiv a
                        }
                (  \sgn m) \zn^{b m}
        \cr& = {- 1 \over N} \sum_{m\divides n} (\sgn m) \zn^{k  n / m} \sumb (- \zn^m
                )^b
        \cr&= - \sum_{
                \scriptstyle m \divides  n
                \atop
                \scriptstyle m \equiv N / 2
                }
                (\sgn m) \zn^{k  n / m}
\cr}
$$
We see that $\gamma_ n = 0$ unless $ {N \over 2} \divides  n $, and we find that
$$
\eqalign{
        \gamma_{ n N / 2}
        & = -\sum_{
                \scriptstyle m\divides n
                \atop
                \scriptstyle m \text{\hskip 3pt odd}
                }
                (\sgn m) \zn^{k  n / m}
        \cr& = - \sum_{
                \scriptstyle m\divides n
                \atop { \scriptstyle m \text{\hskip 3pt odd}
                        \atop
                        \scriptstyle m > 0
                      }
                }
                \zn^{k  n / m} - \zn^{- k  n / m},
}
$$
which is, indeed, the coefficient of $q^{ n/2}$ in $h_k(\tau)$.
\endproof

We make the connection with modular forms by observing that
$$
h_k(\tau) =
G'\bigl(
   \Z \cdot \tpi + \Z \cdot \tpi \tau,
   \tpi \bigl({k \over N} + \half \tau\bigr)\bigr).
\tag 3.6
$$
which together with (3.4) justifies (2.5).

\section 4. Expansion at the cusp $\tau=0$.

In order to test the rationality of $G'$ we will examine the $q$-expansion at a
cusp which happens to be a rational point on the curve $X_1(N)$.
We may, for example, expand
$$
G'(\ztpi+\ztpitau,\tpi{\ell\tau\over N})
$$
in terms of $q^{1/N} = e^{\tpi\tau/N}$, where $\ell$ is any integer relatively
prime to $N$.  Only one such cusp is needed to apply the $q$-expansion
principle, so we could set $\ell=1$, but the computation is no harder if we
refrain from doing that.

Proceeding as in Proposition~3.2, we may write down
the expansions we desire in the following proposition; they don't seem to
follow easily from \cite{14, 2.7.8}.

\proclaim{Proposition 4.1}
\itemitem{(a)}
 $
 G'(\ztpi+\ztpitau,\tpi{\ell\tau\over N})
 = ({\ell \over N} - \half) -
                \sum\limits_{ n=1}^\infty \biggl(
                        \sum\limits_{
                                \scriptstyle m \divides  n
                                \atop
                                \scriptstyle m \equiv \ell
                                }
                             \sgn m \biggr) q^{ n/N}$
\itemitem{(b)}
 $
 G'(\ztpi+\ztpitau,\tpi( \half + {\ell\tau\over N} ))
 = ({\ell \over N} - \half) -
                \sum\limits_{ n=1}^\infty \biggl(
                        \sum\limits_{
                                \scriptstyle m \divides  n
                                \atop
                                \scriptstyle m \equiv \ell
                                }
                             (-1)^{ n/m} \sgn m \biggr) q^{ n/N} $
\endproclaim

\proof
We have
$$
\eqalign{
G'(\ztpi+\ztpitau,\tpi{ \ell\tau \over N} )
& = \ootpi G'(\Z+\Z\tau,\ell\tau/N)
\cr&= \ootpi \sumab \chi^{\Z+\Z\tau}_{\ell\tau/N}\({a \tau + b \over N}\)
                G(\Z+\Z\tau, { a \tau + b \over N } )
\cr&= \ootpi \sumab \zn^{-\ell b} \Gab(\tau),
\cr
}
$$
so applying (3.1) we can write
$$
G'(\ztpi+\ztpitau,\tpi\ell\tau/N)
= \sum^\infty_{ n=0} \beta'_ n q^{ n/N}
$$
where
$$
\beta_\nu' = { 1 \over \tpi} \sumab \zn^{- \ell b} \alpha_ n(a,b)
$$
and
$$
 G'(\ztpi+\ztpitau,\tpi( \half + {\ell\tau\over N} ))
=
 \sum_{ n=0}^\infty \gamma'_ n q^{ n/N},
$$
where
$$
\gamma'_ n = {1
\over \tpi} \sumab (-1)^a \zn^{- \ell b} \alpha_ n(a,b).
$$  We compute
$$
\eqalign{
        \beta_0'
        &= {1 \over \tpi} \sumab \zn^{- \ell b} {1 \over N}
                \,\delta\!\left( { a \over N } \right)
                \lim_{s \rightarrow 1} [
                        \zeta(s,b/N) - \zeta(s,-b/N)
                        ]
\cr
        &= {1 \over \tpi N } \lim_{s \rightarrow 1} \sumb
                        \zn^{- \ell b} [
                        \zeta(s,b/N) - \zeta(s,-b/N)
                        ]
\cr
        &= {-1 \over \tpi N } \lim_{s \rightarrow 1} \sumb
                        (\zn^{\ell b} - \zn^{- \ell b})
                        \zeta(s,b/N)
\cr
        &= {-1 \over \tpi N } \lim_{s \rightarrow 1} \sum_{b=1}^N
                        (\zn^{\ell b} - \zn^{- \ell b})
                        \sum_{n=0}^\infty\left(n + {b \over N}\right)^{-s}
\cr
        &= {-1 \over \pi N } \lim_{s \rightarrow 1} N^s \Im
                \sum_{n=0}^\infty \sum_{b=1}^N \zn^{\ell b} \left(N n + b\right)^{-s}
\cr
        &= {-1 \over \pi} \lim_{s \rightarrow 1} \Im
                \sum_{m=1}^\infty \zn^{\ell m} m^{-s}
\cr
        &= {-1 \over \pi} \Im \sum_{m=1}^\infty \zn^{\ell m} m^{-1}
\cr}
$$
The latter series converges because $\sumb \zn^{\ell b} = 0$,
which tells us that its terms, when taken $N$ at a time, are of the order
of $m^{-2}$.
The last equality used above is justified by the continuity of Dirichlet series
up to the line of convergence, \cite{6, p.~87}, provided we let $s$
approach $1$ from the right through real values.
The last series above is summed in formula (7) of \cite{26, p.~54};
alternatively, we may
apply the continuity of power series up to the circle of
convergence (analogous to the continuity of Dirichlet series up to the line of
convergence), \cite{1, Lehrsatz~IV}, letting $z$ approach $\zn^\ell$ radially from the origin, we obtain
$$
\eqalign{
\beta_0'
&= { -1 \over \pi } \lim_{z \to \zn^\ell} \Im \sum_{m=1}^\infty z ^ m m^{-1}
\cr
&= {  1 \over \pi } \lim_{z \to \zn^\ell} \Im \log (1 - z)
\cr
& = { 1 \over \pi } \arg \left( 1 - \zn^\ell \right)
\cr
& = { \ell \over N } - \half
\cr}
$$

        Now for $ n \geq 1 $ we compute
$$
\eqalign{
        \beta'_ n
        &= { 1 \over \tpi } \sumab \zn^{- \ell b} \biggl[
                - { \tpi \over N } \sum_{
                        \scriptstyle m\divides n
                        \atop
                        \scriptstyle {  n \over m } \equiv a
                        }
                ( \sgn m ) \zn ^ {b m} \biggr]
        \cr&= { - 1 \over N } \sum_{m\divides n} (\sgn m) \sumb \zn^{b(m-\ell)}
        \cr&= - \sum_{
                \scriptstyle m \divides  n
                \atop
                \scriptstyle m \equiv \ell
                } \sgn m
\cr }
$$
This proves (b).  The proof for (c) is similar. \endproof

\demo{\bf Alternate proof of 4.1}
We offer another proof based on formula {\bf H6} of \cite{18,
p.~250}; it uses the Weierstrass $\zeta$-function.  For the purpose of
reconciling the notation in \cite{18, pp.~247--250} with that in
\cite{14, Appendix~C}, we imagine that each $\eta$ in
\cite{18} has been replaced by $-\eta$.  We use the notations
$G_{1,\ell}$, $G^{1,\ell}$, $V_N$, $F$, $h_{a_1,a_2}$, and $H_{a_1,a_2}$ from
\cite{18} without repeating the definitions, and we assume $0 < \ell <
N$.

From \cite{18, {\bf H6}, p.~250} and the definition of $G_{1,\ell}$ we
have
$$
G^{1,\ell}(\tau) = {\ell \over N} - \half
        - \sum_{ j=0}^\infty {
            q^{( j N + \ell) / N} \over 1 - q^{( j N + \ell) / N}
                                }
        + \sum_{ j=1}^\infty {
            q^{( j N - \ell) / N} \over 1 - q^{( j N - \ell) / N}
                                }.
$$
By collecting similar powers of $q$ we obtain
$$
G^{1,\ell}(\tau) = {\ell \over N} - \half
        - \sum_{ n=1}^\infty \biggl(
            \sum_{ \scriptstyle  m\divides n \atop \scriptstyle m \equiv \ell}
                \sgn m
                \biggr) q^{ n/N}.
$$
On the other hand, we may simply trace the definitions in \cite{18}
and \cite{14}.
$$
\eqalign{
      G^{1,\ell}(\tau)
   &= \ootpi h_{\ell/N,0}(\tau)
\cr&= \ootpi H_{\ell/N,0}(\tau,1)
\cr&= \ootpi \(   \zeta(\ell\tau/N,\Z+\Z\cdot\tau)
                + \eta(\ell\tau/N,\Z+\Z\cdot\tau) \)
\cr&= \ootpi \(   \zeta(\ell\tau/N,\Z+\Z\cdot\tau)
                + {1 \over N} \eta(\ell\tau,\Z+\Z\cdot\tau) \)
\cr&= \ootpi A_1(\Z+\Z\cdot\tau,\ell\tau/N)
\cr&= A_1(\ztpi+\ztpitau,\tpi\ell\tau/N)
\cr&= G'(\ztpi+\ztpitau,\tpi\ell\tau/N).
\cr}
$$
Combining the previous two formulas gives (a).

As for (b), we see as above that
$$
G'(\ztpi+\ztpitau,\tpi(\half + {\ell\tau \over N}))
=
\ootpi h_{\ell/N,1/2}(\tau).
$$
Then we apply {\bf H4} of \cite{18}:
$$
\eqalign{
\ootpi h_{\ell/N,1/2}(\tau)
   &= {\ell \over N} + F(q,-q^{\ell/N})
\cr&= {\ell \over N} - \half
        - \sum_{ j=0}^\infty {
                -q^{( j N + \ell)/N} \over 1 + q^{( j N + \ell)/N}
                }
        + \sum_{ j=1}^\infty {
                -q^{( j N - \ell)/N} \over 1 + q^{( j N - \ell)/N}
                }
\cr&= {\ell \over N} - \half
        - \sum_{ n=1}^\infty \biggl(
                \sum_{
                \scriptstyle m \divides  n
                \atop
                \scriptstyle m \equiv \ell
                } (-1)^{n/m} \sgn m
                                \biggr) q^{ n/N}
\cr}
$$
\endproof

{\bf Remark 4.2.}
One may also follow Lang \cite{18,{{\bf H4}}} as in the alternate
proof
above to get the $q$-expansions in (3.3) and (3.6).
One uses
$$
\eqalign{
      A_1(&\ztpi + \ztpitau, \tpi k / N)
\cr&= \ootpi h_{0,k/N}(\tau)
\cr&= F(q,\zn^k)
\cr&= - \half - {\zn^k \over 1-\zn^k}
        - \sum_{ j=1}^\infty{ q^ j\zn^k \over 1 - q^ j\zn^k }
        + \sum_{ j=1}^\infty{ q^ j\zn^{-k} \over 1 - q^ j\zn^{-k} }
\cr&= \half\({\zn^k+1 \over \zn^k-1}\)
        - \sum_{ n=1}^\infty \biggl(
                \sum_{
                \scriptstyle m \divides  n
                \atop
                \scriptstyle m > 0
                }
                \(\zn^{km}-\zn^{-km}\)
                \biggr) q^ n
\cr&= g_k(\tau)
\cr}
$$
and
$$
\eqalign{
      A_1(&\ztpi + \ztpitau, \tpi ( { k \over N} + \half \tau ))
\cr&= \ootpi h_{1/2,k/N}(\tau)
\cr&= \half + F(q,q^{1/2}\zn^k)
\cr&= \half - \half
   - \sum_{ j=0}^\infty{ q^{ j+1/2}\zn^k \over 1 - q^{ j+1/2}\zn^k }
   + \sum_{ j=1}^\infty{ q^{ j-1/2}\zn^{-k} \over 1-q^{ j-1/2}\zn^{-k} }
\cr&= - \sum_{ n=1}^\infty \biggl(
                \sum_{
                \scriptstyle 0 < r \divides  n
                \atop
                \scriptstyle r \text{\hskip 3pt odd}
                } \(
                        \zn^{k n/r}-\zn^{-k n/r}
                        \)
                \biggr) q^{ n/2}
\cr&= h_k(\tau)
\cr}
$$

\section 5. The $q$-expansion principle.

The $q$-expansion principle is part of the classical theory of modular forms.
When the modular curve in question has a model over a subfield $k \subseteq \C$,
a modular function tends to be defined over $k$ if and only if the
coefficients of its $q$-expansion are in $k$.  In this section we state the
purely algebraic portion of the $q$-expansion principle.

Let $K/k$ be a field extension.  In our application, it will be $\C/\Q$.

Suppose $X$ is a curve over $K$, with $x$ a nonsingular point on $X$ defined
over $K$,
and let $q$ be an element of valuation $1$ in the complete local ring $\widehat \O_x$.
We have a $K$-algebra isomorphism $\widehat \O_x \cong K[[q]]$ with the ring $K[[q]]$ of formal
power series in the variable $q$, and an inclusion of fields $K(X)
\hookrightarrow K((q))$ into the field of formal Laurent series.  Pick
generators $f_1,\dots,f_r$ over $K$ for the
field $K(X)$, so that $K(X) = K(f_1,\dots,f_r)$.

\proclaim{Proposition 5.1}
  With the notation above, if $f_1, \dots, f_r \in k((q))$, then
\itemitem{(a)}
        $k((q))$ and $K$ are linearly disjoint over $k$.
\itemitem{(b)}
        $K(X) \cap k((q)) = k(f_1,\dots,f_r)$.
\itemitem{(c)}
        If $X$ has a model $Y$ over $k$ with respect to which the functions
        $f_1,\dots,f_r$ are defined over $k$
        (i.e., $f_i \in k(Y)$),
        then $k(Y) = k(f_1,\dots,f_r)$.
\endproclaim

When we say that $X$ has a model $Y$ over $k$, we mean $Y$ is an irreducible curve
over $k$ with an isomorphism
$g : Y \otimes_k K \isom X$.

\proof
First notice that $k(f_1,\dots,f_r) \subseteq k((q))$.
To prove (a) we proceed as in \cite{23, p.~141}.  Given $a_1,\dots,a_m \in K$
linearly independent over $k$, we suppose that $g_1,\dots,g_m \in k((q))$
satisfy
$\sum a_i g_i = 0$.  We write each $g_i$ as a Laurent series $g_i = \sum_j b_{ij}
q^j$ with coefficients $b_{ij} \in k$.  Then we deduce that $\sum_i \sum_j a_i
b_{ij} q^j = 0$, and thus, for each $j$, we have $\sum_i a_i b_{ij} = 0$, whence
each $b_{ij} = 0$, and thus each $g_i = 0$.

We claim that
when we have intermediate fields $k \subseteq L \subseteq M \subseteq
k((q))$, then $K \cdot L \subseteq K \cdot M$, with equality holding iff $L=M$.
We need only prove that $K\cdot L = K\cdot M$ implies $L=M$.
We know $K$ and $k((q))$ are linearly disjoint subfields of $K((q))$, so it
follows that so are $K$ and $L$, since $L$ is a subfield of $k((q))$.
It follows that $K \otimes L$ maps isomorphically
onto a subring $K*L$ of the field $K\cdot L$ and generates that field,
so $K\cdot L$ is the fraction field of $K*L$.  Similarly, $K\cdot M$ is the
fraction field of $K*M$.  But $M$ is a field extension of $L$, hence a free $L$-module,
a property which is preserved by
tensor product.  Hence $K*M$ is a free $K*L$-module, and since
$K\cdot L = K\cdot M$, we can view $K*M$ as a submodule of the
fraction field of $K*L$.  But any free submodule of the fraction field of a ring has rank
at most $1$, so $K*L = K*M$. Thus
$K \otimes L \isom K \otimes M$, implying $L = M$, proving the claim.

To prove (b), we apply this claim with $L = k(f_1,\dots,f_r)$ and $M = K(f_1,\dots,f_r) \cap k((q))$.

The same fact proves (c). Indeed, the isomorphism $g$ above induces an injective map from
$k(Y) \otimes_k K$ into $K(X)$ whose image generates the field $K(X)$.
This implies that $k(Y)$ and $K$ are linearly disjoint as subfields of $K(X)$.
As $k(f_1,...,f_r)$ is a subfield of $k(Y)$, it is also linearly disjoint from $K$.
Now we are done as before.

\endproof

We refer to \cite{7}, \cite{12}, \cite{13} and
\cite{15} for discussions of the $q$-expansion.

\section 6. The function field of the rational model of $\gp\leftquot H$.

We recall the definitions of the standard examples of modular forms and
functions for $\Gamma = \sltwo$.
For each even integer $\ell \ge 4$
define a homogeneous function $E_\ell : \L \to
\C$ of degree $-\ell$ by the formula $E_\ell(L) =
\sumprime_{w \in L} w^{-\ell}$.
Define $G_2 = 60 \, E_4$ and $G_3 = 140 \, E_6$.  Let $g_2 = (\tpi)^4 G_2 \circ \phi_k$ and
$g_3 = (\tpi)^6 G_3 \circ \phi_k$; these are modular forms of weight $4$ and $6$
respectively.
 Take $\Delta = g_2{}^3 - 27 g_3{}^2$  to be the standard cusp form of weight $12$,
and let $j = g_2{}^3/\Delta$ be the $j$-invariant, a modular function.

It is known that the $q$-expansion of $j$ has rational coefficients; since $j$
is a modular function for the full modular group $\Gamma$, the same is true for
its $q$-expansion at the other cusps.

The Weierstrass $\wp$-function gives a modular form of weight $2$ for $\gp$ as
follows.  We define a function $\wp : \LN \to \C$ by the formula $\wp(L,u) =
\sum_{w\in L}\((w-u)^{-2} - w^{-2}\)$; it is a homogeneous function of degree
$-2$.  We let $p_k = \wp \circ \phi_k$ denote the corresponding modular form.
According to \cite{23, 6.2.1, p.~141} the $q$-expansion of $\wp$ is
$$
\eqalign{
\wp(\Z\cdot\tau + \Z, {}&{ r \tau + s \over N}) = (\tpi)^2 \biggl\{
 {1 \over 12}
- 2 \sum_{n=1}^{\infty} { n q^n \over 1-q^n}
\cr&{}
+ \zn^s q^{r/N} \( 1 - \zn^s q^{r/N}\)^{-2}
+ \sum_{n=1}^{\infty}
     \( \zn^{ns} q^{nr/N} + \zn^{-ns} q^{-nr/N} \) {n q^n \over 1-q^n } \biggr\},
\cr}
$$
for $0 \le r < N$, $(r,s) \in \Z^2$, and $(r,s) \notin N \Z^2$.
Taking $r=k$ and $s=0$ we find that the $q$-expansion of
$$
        \wp(\ztpi + \Z \cdot\tpi \tau,\tpi k \tau/N) ,
$$
which is the $q$-expansion for $p_k(\tau)$ at $\tau = 0$, has rational
coefficients.

We define $f_k = g_2 g_3 p_k / \Delta$; it is a modular function for $\gp$.  It
is known that $(\tpi)^{-4} g_2$, $(\tpi)^{-6} g_3$, and $(\tpi)^{-12}
\Delta$ have
$q$-expansions with rational coefficients (see \cite{23, 2.2}).
It follows that the $q$-expansion at the cusp $\tau=0$ of the modular function
$f_k$ has rational coefficients.

\proclaim{Proposition 6.1}
 The field $\C(j,f_k)$ is the field of all modular functions for $\gp$.
\endproclaim

\proof
Let $M$ be the field of all modular functions for the group
$$
\Gamma(N) = \biggl\{ \pmatrix a&b\cr c&d\endpmatrix \in \sltwo
 \biggm|
    \pmatrix a&b\cr c&d\endpmatrix \equiv \pmatrix 1&0\cr 0&1\endpmatrix \mod N
 \biggr\}.
$$
(This field $M$, called the modular function field of level $N$, is the function field of the modular curve $X(N)$ over $\C$.)
We know that $\C(j)$ is the field of all modular functions for the group $\Gamma = \sltwo$,
that $M$ is a Galois extension of $\C(j)$,
and that $\Gal(M/\C(j)) = \Gamma/\Gamma(N)\cdot\{\pm 1\} \simeq$ SL$_2(\Z/N)/\{\pm 1\}$.

Let $L$ be the field of all modular functions for $\gp$.  We know that
$K = \C(j,f_k) \subseteq L$, and that $\Gal(M/L) = \gp \cdot \{ \pm1 \} / \Gamma(N) \cdot \{ \pm1 \}$.
To determine the group
corresponding to the intermediate field $K$ we consider an element $\gamma \in
\Gamma$ which leaves every element of $K$ fixed.  Then, as in \cite{23,
6.1-A}, we see that the equation $f_k \circ \gamma = f_k$ implies that
$$
\pmatrix0&\frac k N\endpmatrix \cdot \gamma \equiv \pm \pmatrix0 & \frac k
N\endpmatrix
\mod{\Z^2},
$$
and thus that
$$
\gamma \equiv \pmatrix *&*\\0&\epsilon\endpmatrix \mod N,
$$
where $\epsilon = \pm 1$.  The matrix $\gamma' = \epsilon \gamma$ has
determinant $1$, so
$$
\gamma' \equiv \pmatrix1&*\\0&1\endpmatrix \mod N,
$$
i.e., $\gamma' \in \gp$, showing that $\gamma \in \gp \cdot \{\pm1\}$.  Thus $L = K$.
\endproof

Applying (5.1) yields the following corollary.

\proclaim{Corollary 6.2}
 The field $\Q(j,f_k)$ is the field of modular functions for
$\gp$ whose $q$-expansions at the cusp $\tau=0$ have rational coefficients, and
is the function field of the canonical rational model of $X_1(N)$.
\endproclaim

\bigskip


\Refs

\newfam\cyrfam

\font\eightcyr=wncyr8
\input cyracc.def
\def\cyr{\eightcyr\cyracc}



\loadmsbm
\loadmsam
\font\eightmsb=msbm8
\textfont\msbfam=\eightmsb

\input amsppt.sty 
\raggedright

\catcode`\@=11
\def\logo@{}
\catcode`\@=\active
\tolerance=10000

\def\num#1{}
\def\numref#1{}

\def\callnumber#1{\finalinfo [#1]}

\hyphenation{pre-print}



\def\birkhaeuser{\publ Birkh\"auser \publaddr Boston, Basel, Berlin}

\def\ams{\publ American Mathematical Society \publaddr Providence,
	Rhode Island}
\def\springer{\publ Springer \publaddr Berlin, Heidelberg, New York}


\def\lnm{Lecture Notes in Mathematics}
\def\gradtext#1.{\bookinfo Graduate Texts in Mathematics \vol #1 \springer}

\def\PM{\bookinfo Progress in Mathematics \birkhaeuser}

\def\cnm{Contemporary Mathematics}

\def\CNM{\ams\bookinfo\cnm\ }

\def\grundlehren#1.{\springer
	\bookinfo Grundlehren der mathematischen Wissenschaften \vol #1}
\def\sga#1:#2.{\bookinfo S\'eminaire de G\'eom\'etrie Alg\'ebrique
	du Bois-Ma\-rie #2, {\bf SGA #1}}
\def\sgalnm#1:#2:#3.{\bookinfo S\'eminaire de G\'eom\'etrie Alg\'ebrique
	du Bois-Ma\-rie #2, {\bf SGA #1}, \lnm \vol #3 \springer}


\def\bOulder#1{\CNM 55
	\inbook Applications of Algebraic $K$-theory to
	Algebraic Geometry and Number Theory, Part #1,
	{\rm Proceedings of a Summer Research Conference held
		June 12-18, 1983, in Boulder, Colorado}
	}
\def\BoulderI{\bOulder{I}}

\def\AriAlgGeo{\PM \vol 89 \inbook Arithmetic Algebraic Geometry (proceedings
	of a conference on Texel Island, April, 1989)
	\eds G. van der Geer, F. Oort, and J. Steenbrink \yr 1991
	\callnumber{516.35 Ar462}}





\def\invm{\jour Inventiones Mathematicae}



\ref \no 1 \by N. H. Abel
        \paper Untersuchungen \"uber die Reihe:
        $1 + {m \over 1} \cdot x
           + { m \cdot (m-1) \over 1 \cdot 2} \cdot x^2
           + {m \cdot (m-1) \cdot (m-2) \over 1 \cdot 2 \cdot 3 } \cdot x^3
           + \cdots \text{u.s.w.}$
        \jour Journal f\"ur Math.{}
        \yr 1826 \pages 311--339
        \endref

\ref \no {\relax 2} \by A. Beilinson
        \paper Higher regulators of modular curves
        \book Applications of Algebraic K-theory to
                Algebraic Geometry and Number Theory, Part I
        \bookinfo Contemporary Mathematics \vol 55, part I \yr 1986
        \publ American Math. Soc. \pages 1--34
        \endref

\ref \no {\relax 3} \by Spencer J. Bloch
        \book Higher regulators, algebraic $K$-theory, and zeta functions of elliptic curves
        \publ Amer.\ Math.\ Soc.{}
        \publaddr Providence
        \yr 2000
        \endref

\ref \no {\relax 4} \by S. Bloch and D. Grayson
        \paper $K_2$ and L-functions of elliptic curves :\ Computer Calculations
        \book Applications of Algebraic K-theory to
                Algebraic Geometry and Number Theory, Part I
        \bookinfo  Contemporary Mathematics \vol 55
        \publ Amer.\ Math.\ Soc.{}
        \publaddr Providence
        \yr 1983
        \pages 79--88
        \endref

\ref \no 5 \by C. Breuil, B. Conrad, F. Diamond and R. Taylor
        \paper On the modularity of elliptic curves over $\Q$: wild $3$-adic exercises
        \jour J. Amer. Math. Soc.
        \yr 2001
        \pages 843--939
        \vol 14
        \endref

\ref \no 6 \by E. Cahen
        \paper Sur la fonction $\zeta(s)$ de Riemann et sur des fonctions
                analogues
        \jour Ann. de l'\'Ecole Normale
        \yr 1894 \pages 75--164
        \vol 11 (ser.~3)
        \endref

\ref \no {\relax 7} \by P. Deligne, N. Rapoport
        \book Modular functions of One Variable II
        \bookinfo Lecture Notes in Mathematics
        \vol 349
        \publ Springer-Verlag
        \publaddr Berlin, Heidelberg, New York
        \endref

\ref \no 8
        \by C. Deninger
        \paper Higher Regulators and Hecke $L$-series
        of imaginary quadratic fields: I \invm \vol 96 \yr 1989 \pages 1-69
        \endref

\ref \no 9 \by V. G. Drinfeld
        \paper Two theorems on modular curves
        \jour Funktsionalnii anal. i evo Prilozhenie
        \vol 2
        \yr 1973
        \pages 83--84
        \endref

\ref \no {\relax 10} \by E. Hecke
        \paper Darstellung von Klassenzahlen als Perioden von Integralen
                3.~Gattung aus dem Gebiet der elliptischen Modulfunctionen
        \jour Abhandlungen aus dem Mathematischen Seminar der Hamburgische Universit\"at
        \vol 4
        \yr 1925
        \pages 211--223
        \endref

\ref \no {\relax 11} \by E. Hecke
        \paper Theorie der Eisensteinschen Reihen h\"oherer Stufe und ihre
                Anwendung auf Funktionentheorie und Arithmetik
        \jour Abhandlungen aus dem Mathematischen Seminar der Hamburgische Universit\"at
        \vol 5
        \yr 1927
        \pages 199--224
        \endref

\ref \no {\relax 12} \by N.~Katz
        \paper $p$-adic interpolation of real analytic Eisenstein series
        \jour Annals of Mathematics
        \yr 1976 \vol 104 \page 459
        \endref

\ref \no 13 \by N.~Katz
        \paper $p$-adic $L$-functions for CM-fields
        \jour Inventiones Mathematicae
        \yr 1978 \vol 49
        \endref

\ref \no {\relax 14} \by N.~Katz
        \paper The Eisenstein measure and p-adic interpolation
        \jour American Journal of Mathematics
        \yr 1977 \vol 99 \pages 238--311
        \endref

\ref \no {\relax 15} \by N.~Katz and B.~Mazur
        \paper Arithmetic moduli of elliptic curves
        \jour Annals of Math. Studies
        \vol 108
        \yr 1985
        \pages 1--514
        \endref

\ref \no {\relax 16} \by V.~Kolyvagin
        \paper Finiteness of $E(\Q)$ and ${\cyr Sh}(E,\Q)$ for a subclass of Weil curves
        \jour  Izv. Akad. Nauk SSSR Ser. Mat.
        \vol 52
        \yr 1988
        \pages 522--540, 670--671
        \endref

\ref \no {\relax 17} \by D.~Kubert
        \paper Universal bounds on the torsion of elliptic curves
        \jour Proc.\ London Math.\ Soc.
        \vol 33
        \yr 1976
        \pages 193--237
        \endref

\ref \no {\relax 18} \by S.~Lang
        \book Introduction to Modular Forms
        \bookinfo Grundlehren der mathematischen Wissenschaften
        \publ Spring\-er Ver\-lag
        \publaddr Berlin, Heidelberg, New York
        \yr 1976
        \endref

\ref \no {\relax 19} \by B.~Mazur
        \paper Rational Isogenies of prime degree
        \jour Inventiones Mathematicae
        \yr 1978
        \vol 44
        \pages 129--162
        \endref

\ref \no 20
        \AriAlgGeo \by J.-F. Mestre and N. Schappacher \paper S\'eries de
        Kronecker et fonctions $L$ des puissances sym\'etriques de
        courbes elliptiques sur $\Bbb Q$ \pages 209-246 \endref

\ref \no 21
        \by D. Ramakrishnan \paper Analogs of the Bloch-Wigner function for
        higher polylogarithms \BoulderI \pages 371-376 \endref

\ref \no {\relax 22} \by B.~Schoeneberg
        \paper Elliptic Modular Functions
        \publ Springer
        \publaddr Berlin, Heidelberg, New York
        \yr 1974
        \endref

\ref \no {\relax 23} \by G.~Shimura
        \paper Introduction to the arithmetic theory of
                automorphic functions
        \publ Iwanami Shoten, Publishers, and Princeton University Press
        \yr 1971
        \endref

\ref \no {\relax 24} \by J. Silverman
        \book The Arithmetic of Elliptic Curves
        \publ Springer-Verlag
        \publaddr New York, Berlin, Heidelberg, Tokyo
        \bookinfo Graduate Texts in Mathematics {\bf 106}
        \yr 1986
        \endref

\ref \no {\relax 25} \by R. Taylor and A. Wiles
        \paper Ring-theoretic properties of certain Hecke algebras
        \jour Ann. of Math. (2)
        \vol 141
        \yr 1995
        \pages 553--572
        \endref

\ref \no {\relax 26} \by A. Weil
        \book Elliptic Functions according to Eisenstein and Kronecker
        \bookinfo Ergebnisse series {\relax 88}
        \publ Springer
        \publaddr Berlin, Heidelberg, New York
        \yr 1976
        \endref

\ref \no {\relax 27} \by A. Wiles
        \paper Modular elliptic curves and Fermat's last theorem
        \jour Ann. of Math. (2)
        \vol 141
        \yr 1995
        \pages 443--551
        \endref

\ref \no 28
        \by D. Zagier \paper Polylogarithms, Dedekind Zeta functions and the algebraic $K$-theory of fields
        \jour Progress in Math. \vol 89 \yr 1991 \pages 391--430
        \endref

\endRefs

\bye